# Solving Sangaku With Traditional Techniques

## Rosalie Hosking

**Abstract.**

Between 17th and 19th centuries, mathematically orientated votive tablets appeared in Shinto shrines and Buddhist temples all over Japan. Known as *sangaku*, they contained problems of a largely geometrical nature. In the 17th century, the Japanese mathematician Seki Takakazu developed a form of algebra known as *tenzan jutsu*. I compare one mathematical problem from the 1810 Japanese text *Sanpō Tenzan Shinan* solved using *tenzan jutsu* to a similar problem found on the Kijimadaira Tenman-gū shrine *sangaku* to show how *sangaku* problems can be solved using the traditional Japanese methods.

## §1. Introduction

During the Edo Period (1603-1867 CE), the serene temples and shrines of Japan began to display amongst their artifacts a perhaps unexpected treasure - votive mathematical tablets. These tablets were known as *sangaku* 算額. The oldest surviving tablet dates back to 1683, and historians estimate that thousands were produced during the Edo and early Meiji (1868 - 1912) periods [1, p. 9].

A typical *sangaku* is constructed of wood - preferably paulownia, which did not bend or warp easily - and colourfully presented using bright pigment paints. They ranged in size, some being 20 cm in length and others over 1 metre. The traditional Japanese text they contain is read from right to left and top to bottom. Usually they are written in the *kanbun* 漢文 language, a form of classical Chinese with Japanese readings. Tablets may be divided into several sections, each dealing with a geometrical problem. These problems challenge the reader to find some length, area, or diameter relating to the diagram in terms of other magnitudes. Though instructions for calculation are provided, they rarely offer any clue as to how one might actually solve the problem, for they are formulaic, and rely on simplified arithmetic procedures, giving little reference to the geometrical shapes and relations in question. Because of this, it is difficult for historians to understand and recreate the original methods used by *sangaku* mathematicians to solve these problems.

However, through examining textbooks of the Edo period in conjunction with *sangaku*, clues can be found to understand how these problems may have been solved using traditional methods available at the time. To illustrate this idea I examine one problem from a *sangaku* dedicated to the Kijimadaira Tenman-gū shrine 木島平天満宮 and a similar problem from the Edo period text *Sanpō Tenzan Shinan* 算法点竄指南 *Instruction in Tenzan Mathematics*. I transcribe and translate these two problems into English, and show how the textbook example can be used as a traditional solution for the *sangaku*, allowing us to see how *sangaku* may have been solved by mathematicians of the Edo and Meiji periods.

## §2. The Kijimadaira Tenman-gū Sangaku and Sanpō Tenzan Shinan

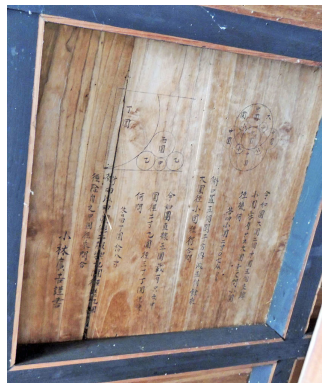

Fig. 1. The Kijimadaira Tenman-gū sangaku

The Kijimadaira Tenman-gū 木島平天満宮 shrine is dedicated to Tenjin 天神, the Japanese deity of education and science. It is located on a wooded hillside on the outskirts of Kijimadaira 木島平, a small rural village in Nagano prefecture. The tablet which concerns us is one of three in this shrine, all which are believed to have been dedicated around 1888 [2, p. 5][1]. The tablet measures 45 cm by 45 cm and forms one panel of the chequered roof of the shrine. The author was Kobayashi Hirokichi 小林廣吉 (?), a local artisan, and his tablet is presented in the typical *sangaku* style.

The problem examined on this *sangaku* is very similar to one found in the 1810 mathematical treatise the *Sanpō Tenzan Shinan* by Ōhara Toshiaki 大原利明 (?-1828). The *Sanpō Tenzan Shinan* is a three volume collection detailing how to solve mathematical problems using the Edo period symbolic manipulation method known as *tenzan jutsu* 点竄術. By examining how this problem from the *Sanpō Tenzan Shinan* is solved, we can ascertain one possible original method for solving the Kijimadaira Tenman-gū *sangaku* using traditional techniques.

§**3. Tenzan Jutsu**

*Tenzan jutsu* is a symbolic manipulation technique for solving equations with more than one unknown. It developed out of the *bōsho hō* 傍書法 - often referred to as the method of side-writing - authored by the famed Edo period mathematician Seki Takakazu 関孝和 (? - 1708). *Tenzan jutsu* solves problems with unknowns through a symbolic style of reasoning. It uses vertical lines to represent numerical values and symbols for unknowns. For instance | 甲 expresses 1 甲, || 甲 expresses 2 甲, and so on. To represent division a symbol is placed to the right such that 乙 | 甲 expresses 甲 ÷ 乙. Other operators are indicated with certain characters. Unknowns were represented using characters from the Chinese celestial stems - 甲 *kō*, 乙 *otsu*, 丙 *hei*, 丁 *tei*, 戊 *bo*, 己 *ki*, 庚 *kō*, 辛 *shin*, 壬 *jin*, and 癸 *ki* - and the twelve signs of the Chinese zodiac - 子 *ne* (rat), 丑 *ushi* (ox), 寅 *tora* (tiger), 兎/卯 *u* (rabbit), 辰 *tatsu* (dragon), 巳 *mi* (snake), 午 *uma* (horse), 未 *hitsuji* (goat), 申 *sara* (monkey), 酉 *tori* (rooster), 戌 *inu* (dog), and 亥 *i* (boar). The main elements of a diagram (circles, squares, cylinders, etc) were labelled with celestial stems, while auxiliary lines added to diagrams were labelled with signs of the zodiac.

In *tenzan jutsu*, some terms are also labelled with the katakana letters *i* イ, *ro* ロ, *ha* ハ, *ni* ニ, *ho* ホ, *he* ヘ, *to* ト, *chi* チ, and *ri* リ.

---

[1]This dedication date is only given on one tablet, but their similar size, condition, construction, and installation by a professional carpenter has led members of the Kijimadaira community to believe they were all dedicated at the same time.

These are alphabet characters derived from the Heian period (794-1170) poem known as the *iroha*. The poem is a pangram similar to *The quick brown fox jumps over the lazy dog*, which contains every classical Japanese alphabet character exactly once[2]. The labels seem to be applied to terms before and after factorisation to indicate which like terms are to be combined.

### §4. The Kijimadaira Tenman-gū Sangaku

The particular problem from the Kijimadaira Tenman-gū tablet we will examine deals with finding the diameter of a circle within a configuration of overlapping circles. The problem is transcribed and translated below:

PROBLEM DIAGRAM AND TEXT

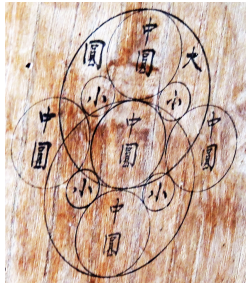

Fig. 2. Kijimadaira Tenman-gū sangaku problem

The $\cdots$ in the answer indicates an approximation has been presented. In modern notation, the technique proposed to solve this problem is $(\sqrt{5} - 5) \times 大$ and the answer is 2.07 (approx). Also the author writes 小径得合問 instead of the more common expression 得小径合問.

### §5. Sanpō Tenzan Shinan Problem

In the *Sanpō Tenzan Shinan*, there is a similar problem to that found on the Kijimadaira Tenman-gū 木島平天満宮 tablet. In this book, we are asked to find the value of the small circles (小) given the diameter

---

[2]The characters *we* ゑ and *wi* ゐ are no longer used in the Japanese language.

| | |
|---|---|
| 今如圖大圓二個中圓五個交罅小圓容有只云大圓壹尺問小圓径幾何 | As in the diagram, there are two large circles (大) and five medium circles (中). In the space between where these circles intersect are small circles (小). Say the diameter of the large circles (大) is 1 *shaku*. Problem - what is the diameter of the small circles (小)? |
| 答曰　小圓二寸0七分有奇 | Answer: The diameter of the small circles (小) is 2 *sun* 0 7 *bu* ⋯. |
| 術曰　置五個開平方內減五個餘乗大圓径小径得合問 | Technique: Put 5 *ko* and take the square root. Inside subtract 5 *ko*. Then multiply by the diameter of the large circles (大). Obtain the diameter of the small circles (小) as required. |

of the large circles (大). However the values in the text are one-tenth of those found on the *sangaku*. In this case the diameter of the large circles is 1 *sun* 寸 rather than 1 *shaku* 尺, where 1 *shaku* = 10 *sun*, and the numerical value for the answer is also one-tenth of that on the *sangaku*, being 2 *bu* 07 *mo*. The formula in the *Sanpō Tenzan Shinan* also reads slightly differently, stating "Put 5 *bu* and take the square root. Inside subtract 5 *bu*. Then multiply by the diameter of the large circles (大). Obtain the diameter of the small circles (小) as required". Here *bu* is used rather than *ko*. The formula can be expressed as $(\sqrt{5bu} - 5bu) \times$ 大 or $(\sqrt{0.5} - 0.5) \times$ 大, with the tablet reading $(\sqrt{5ko} - 5ko) \times$ 大 or $(\sqrt{5} - 5) \times$ 大.

The diagram and solution using *tenzan jutsu* provided by Ōhara in the *Sanpō Tenzan Shinan* are given below in three columns. The first column provides a transcription of the original text, the second a transliteration, and the third a modern expression.

PROBLEM DIAGRAM AND TEXT

[Where 大 = $a$, 中 = $b$, and 小 = $c$.]

| | | |
|---|---|---|
| 解曰置大半之 | | Solution: Put 大 [$a$] and halve it. |
| 二\|大<br><br>中径 | 2\|$a$<br><br>$b$ | $\dfrac{a}{2} = b$ |

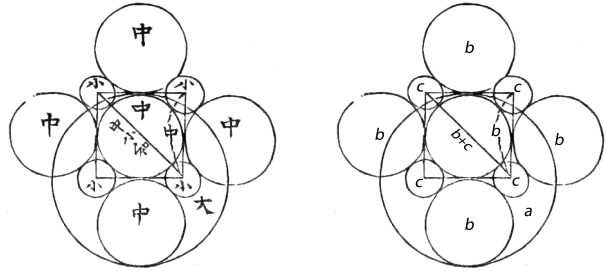

Fig. 3. Left: Diagram from the *Sanpō Tenzan Shinan*. Right: Modern transliteration.

乘方斜率          Multiply by the square diagonal ratio[3]

$$\left. \begin{array}{c|c} 二 & 大\dfrac{二}{商} \\ \hline 中小径和 \end{array} \right. \qquad \left. \begin{array}{c|c} 2 & a\ 2 \\ & \sqrt{} \\ \hline b+c \end{array} \right. \qquad \dfrac{a\sqrt{2}}{2} = b+c$$

乘除省之          Multiplication and division. Eliminate[4]

$$\left. \begin{array}{c|c} \dfrac{二}{商} & 大 \\ \hline 中小径和 \end{array} \right. \qquad \left. \begin{array}{c|c} 2 & a \\ \sqrt{} & \\ \hline b+c \end{array} \right. \qquad \dfrac{a}{\sqrt{2}} = b+c$$

寄左 ◯ 以中小和相消求          Move left. Combine 中 [$b$] and 小 [$c$]. Cancel to find[5]

---

[3]A square is formed in Figure 3 inscribing the middle circle 中. Each side of the square is the length of the diameter 中. Each corner of the square touches the middle of the four circles 小. In this step, the square-diagonal ratio $\sqrt{2}$ is multiplied by the diameter of the circle 中, which is the side length of the square. This produces the length of the diagonal.

[4]The denominator of 2 is converted to $\sqrt{2} \times \sqrt{2}$. Then one instance of $\sqrt{2}$ is cancelled out to produce $\dfrac{a}{\sqrt{2}}$.

[5]The 'cancel' may refer to a cancellation of 小 [$c$] from the right hand side. The character ◯ marks the end of a paragraph.

| | | |
|---|---|---|
| イ 二商 \| 大 | 2 \| a √ | $\dfrac{a}{\sqrt{2}}$ |
| ロ 二 \| 大 | 2 \| a | $-\dfrac{a}{2}$ |
| \| 小 | \| c | $-c = 0$ |
| 合矩 變換得 | 0 Convert to obtain | |
| イ 五大 分商 | 5 a bu √ | $a(\sqrt{0.5})$ |
| 五大 分 | 5 a bu | $-a(0.5)$ |
| 小 | c | $-c = 0$ |
| 合矩 | 0 | |

於是求得小径式    Request an equation satisfied by 小 [c]

| | | |
|---|---|---|
| 五大 分商 | 5 a bu √ | $a(\sqrt{0.5})$ |
| 五大 分 | 5 a bu | $-a(0.5) = c$ |
| 得小径式 | *Obtain formula c* | |

依施答術則如左    By means of the absolute answer technique [the answer is] shown above.

## §6. Applying the *Sanpō Tenzan Shinan* working to the Nagano Tenman-gū Sangaku

The formula provided on the tablet in modern notation reads

(1) $$a(\sqrt{0.5}) - a(0.5) = c$$

In the technique section of the *Sanpo Tenzan Shinan* problem, this is expressed by the equivalent $(\sqrt{0.5} - 0.5) \times a$. Using the original terms provided in the text, this reads $(\sqrt{5bu} - 5bu) \times 大$. However, the tablet problem appears to give the formula $(\sqrt{5} - 5) \times a$, or in the original text $(\sqrt{5ko} - 5ko) \times 大$, which does not produce the value of 2.07 given on the tablet. It is curious that the author has done this, for when checking his calculations it should have become immediately clear that an error occurs when using *ko* instead of *bu*.

One explanation could be that the formula was supposed to be calculated on the *soroban*, meaning can be considered vague and ambiguous but not necessarily incorrect. The formula tells us to put down the number 5, which we then square, subtract 5 from, and multiply by 10. On the *soroban* one row represents ones, the next tens, the next hundreds, and so on. This means 5 has the potential to be expressed 0.05, 0.5, 5, 50, or 500. The user can also choose which column they wish to start from (i.e. use as the ones) or put down a value in any row to treat it as abstract. It is possible that due to this vagueness the author has used different units but obtained the same formula and result.

By examining the *Sanpo Tenzan Shinan*, we can see how a traditional Japanese mathematician could have approached this *sangaku* problem. The solution essentially relies on a recognition of the square and diagonal formed in Figure 3. From this it can be determined that $\frac{a}{2} = b$ and $\frac{a}{\sqrt{2}} = b + c$. By substituting $\frac{a}{2}$ for $b$ the value of $c$ can then be obtained.

It can also be seen from this example that while *sangaku* provide a simple formula, the working required to obtain that formula using traditional methods is not straightforward or easily visible from the tablet. While it cannot be said that *sangaku* authors approached problems the same way as authors of mathematical books such as Ōhara, there is some evidence that the mathematicians of these books produced *sangaku*. For example, the first volume of another Edo period text the *Shinpeki Sanpō* 神壁算法 *Mathematics of God's wall* exclusively deals with problems displayed on tablets dedicated by members of the Seki Takakazu 関孝和 school of mathematics[4, p. 142]. This indicates that *sangaku* did not constitute a separate mathematical tradition, but were a branch of the wider Edo period tradition which employed approaches and methods as found in the *Sanpō Tenzan Shinan*. So it can be seen from this example that by examining *sangaku* in conjunction with texts like the *Sanpō Tenzan Shinan*, more clarity can be found on how to solve *sangaku* problems traditionally with Edo and early Meiji methods.

*School of Mathematics and Statistics, University of Canterbury, Christchurch, New Zealand*
*E-mail address*: `rosalie.hosking@gmail.com`